\documentclass[12pt,reqno]{amsart}

\usepackage{amsmath,amsthm,amsfonts,amssymb,mathtools}
\usepackage{mathrsfs}
\usepackage{graphicx,xcolor}
\usepackage{enumitem}
\usepackage[bookmarksnumbered,colorlinks,plainpages]{hyperref}

\newcommand{\R}{\mathbb R}
\newcommand{\C}{\mathbb C}

\newcommand{\Cp}{\mathbb C^+}
\newcommand{\Cm}{\mathbb C^-}

\newcommand{\Pp}{\mathcal P}

\renewcommand{\d}{\,\mathrm{d}}

\newcommand{\weak}{\Rightarrow}

\newcommand{\Law}{\operatorname{Law}}

\newcommand{\kcx}{\prec_{k\text{-}\mathrm{cx}}}

\newtheorem{theorem}{Theorem}[section]
\newtheorem{lemma}[theorem]{Lemma}
\newtheorem{proposition}[theorem]{Proposition}

\theoremstyle{definition}

\theoremstyle{remark}

\numberwithin{equation}{section}

\begin{document}
\setcounter{page}{1}

% Journal header: fill in after acceptance.
%\color{darkgray}{
%\noindent
%{\small Journal Name}\hfill
%{\small ISSN: XXXX-XXXX}\\
%{\small Vol.~XX (20XX) XX--XX}\hfill
%{\small DOI: to be assigned}}
%\color{black}

\title[Higher-order convex orders under free multiplicative convolution]
{Higher-order convex orders under free multiplicative convolution}

\author[Cheng-Xuan Soong]{Cheng-Xuan Soong}

\address{
National Center for Theoretical Sciences, Mathematics Division,
National Taiwan University, Taipei 106, Taiwan
}

\email{cxsoong@ncts.ntu.edu.tw}

\date{\today}

\begin{abstract}
We prove that free multiplicative convolution with a probability
measure on $[0,\infty)$ preserves the $k$-convex order for every
integer $k\ge2$, assuming finite $(k-1)$st absolute moments for
the two measures being compared and for the convolving measure.
The compared measures may have arbitrary support in $\R$.

%\medskip
%\noindent\textit{2020 Mathematics Subject Classification.}
%Primary 46L54; Secondary 60E15.

\end{abstract}

\maketitle

%%%%%%%%%%%%%%%%%%%%%%%%%%%%%%%%%%%%%%%%%%%%%%%%%%%%%%%%%%%%%%%%%%%%%%
\section{Introduction}\label{sec:intro}
%%%%%%%%%%%%%%%%%%%%%%%%%%%%%%%%%%%%%%%%%%%%%%%%%%%%%%%%%%%%%%%%%%%%%%

Let $\Pp_k(\R)$ denote the set of Borel probability measures on $\R$
with finite $k$th absolute moment, and $\Pp_k(\R_{\ge0})$ the subset
of those supported on $[0,\infty)$; we write $\Pp_0$ when no moment condition is imposed. If $X=X^*$ and $Y\ge0$ are freely independent with distributions $\mu$ and $\nu$, respectively, their free multiplicative convolution, introduced by Voiculescu \cite{Voi87}, is
$\mu\boxtimes\nu=\Law(Y^{1/2}XY^{1/2})$. This operation is defined for
all $\mu\in\Pp_0(\R)$ and $\nu\in\Pp_0(\R_{\ge0})$, including
unbounded supports, in which case $X$ and $Y$ are affiliated
operators; see \cite{BV93,AHK26}.

Higher-order convex orders compare measures by functions with
nonnegative higher-order divided differences; we follow the
formulation of Denuit, Lef\`evre and Shaked \cite{DLS98} (see
Subsection~\ref{subs:order}). For a smooth function, $k$-convexity is
equivalent to nonnegativity of its $k$th derivative, and $k=1,2$ give
the usual stochastic and convex orders. Classical multiplication by
an independent nonnegative random variable preserves the $k$-convex
order. Indeed, the law of the product is a mixture of dilations of
the law of the first factor, and $x\mapsto f(yx)$ is $k$-convex for
every $y\ge0$ whenever $f$ is. In the free setting,
$\mu\boxtimes\nu$ is not, in general, a mixture of dilations of
$\mu$, so this reduction is unavailable. The corresponding
preservation property for free additive convolution was proved in
\cite{Soo26}, whose $k=4$ case answers a question of Hein\"avaara
\cite[Question~5.8]{Hei26}. The purpose of this paper is to show that
the preservation property also holds for free multiplicative
convolution.

\begin{theorem}\label{thm:main}
Let $k\ge2$ be an integer, let $\mu_1,\mu_2\in\Pp_{k-1}(\R)$, and let
$\nu\in\Pp_{k-1}(\R_{\ge0})$. If $\mu_1\kcx\mu_2$, then
\[
    \mu_1\boxtimes\nu\kcx\mu_2\boxtimes\nu .
\]
\end{theorem}

For $k=2$, only a finite first moment is required of $\nu$,
whereas the additive result in \cite[Theorem~1.1]{Soo26}
assumes a finite second moment of the convolving measure.
The support condition $\nu\in\Pp_0(\R_{\ge0})$ is part of
the definition of $\mu\boxtimes\nu$ for $\mu\in\Pp_0(\R)$.
It is also necessary for a preservation statement covering
all orders, even classically: multiplication by $-1$ reverses
the $k$-convex order for odd $k$, since $x\mapsto-f(-x)$ is
$k$-convex whenever $f$ is $k$-convex.

The case $k=1$ requires no moment assumptions and follows directly
from \cite[Proposition~4.15]{BV93}. For compactly supported measures,
the case $k=2$ also admits an operator-theoretic proof, combining
Strassen's martingale coupling theorem \cite{Strassen65},
a free-product construction (cf.~\cite{Boca91}), and spectral
majorization \cite{Hiai87}. This argument does not extend directly
to $k\ge3$, since the matching second moments force any martingale
coupling to be trivial. To treat all $k\ge2$ uniformly, we adapt
the analytic approach of \cite{Soo26}, using the multiplicative
subordination theorem of \cite{AHK26}.

The paper is organized as follows.
Section~\ref{sec:prelim} collects preliminary results on
higher-order potentials, uniqueness of transforms, and
multiplicative subordination.
Section~\ref{sec:variation} derives the variation formula along
a segment of probability measures.
Section~\ref{sec:kernels} constructs positive kernels and
establishes an integrated bound depending only on the first
moment of $\nu$.
Section~\ref{sec:proof} completes the proof of
Theorem~\ref{thm:main}.

\section{Preliminaries}\label{sec:prelim}

Throughout the remainder of the paper, $k\ge2$ is an integer.
We write $m_j(\mu)=\int_\R x^j\d\mu(x)$ and $m_\nu=m_1(\nu)$
whenever these moments are finite. Let $\Cp$ and $\Cm$ denote
the upper and lower half-planes. For $\mu\in\Pp_0(\R)$, its
Cauchy transform and reciprocal Cauchy transform are
\[
    G_\mu(z)=\int_\R\frac{\d\mu(x)}{z-x},\quad
    F_\mu(z)=\frac1{G_\mu(z)},\quad z\in\C\setminus\R.
\]
We use the same notation $G_\rho$ for finite signed Borel
measures $\rho$ on $\R$. We write $C_c(\R)$ for the continuous
functions with compact support and $C_0(\R)$ for the continuous
functions vanishing at both infinities.

\subsection{Higher-order convex orders and potentials}
\label{subs:order}

Following \cite{DLS98}, a function $f:\R\to\R$ is
\emph{$k$-convex} if
\[
    f[x_0,\ldots,x_k]\ge0
    \quad\text{whenever }x_0<\cdots<x_k,
\]
where $f[x_0,\ldots,x_k]$ denotes the divided difference of
order $k$. For $f\in C^k(\R)$, this is equivalent to
$f^{(k)}\ge0$.

For $\mu_1,\mu_2\in\Pp_{k-1}(\R)$, we write
$\mu_1\kcx\mu_2$ if
\[
    \int_\R f(x)\d\mu_1(x)
    \le
    \int_\R f(x)\d\mu_2(x)
\]
for every $k$-convex function
$f\in L^1(\mu_1)\cap L^1(\mu_2)$.
For $\mu\in\Pp_{k-1}(\R)$, define
\[
    U_k^\mu(t)=\frac1{(k-1)!}
    \int_\R(x-t)_+^{k-1}\d\mu(x),\quad t\in\R,
\]
where $a_+=\max\{a,0\}$ denotes the positive part of $a$.
The moment assumption makes these potentials finite and continuous.

\begin{proposition}[{\cite[Theorem~3.2]{DLS98}}]
\label{prop:DLS}
Let $\mu_1,\mu_2\in\Pp_{k-1}(\R)$. Then $\mu_1\kcx\mu_2$
if and only if $m_j(\mu_1)=m_j(\mu_2)$ for $1\le j\le k-1$
and $U_k^{\mu_1}\le U_k^{\mu_2}$ on $\R$.
\end{proposition}

We will use the following properties of potential differences
and a uniqueness result for higher-order Cauchy transforms.

\begin{lemma}[{\cite[Lemmas~3.3--3.4]{Soo26}}]
\label{lem:potentials}
Let $\mu_1,\mu_2\in\Pp_{k-1}(\R)$ have matching moments through
order $k-1$. Then $U_k^{\mu_2}-U_k^{\mu_1}\in C_0(\R)$ and
\[
    G_{\mu_2}(z)-G_{\mu_1}(z)
    =k!\int_\R
      \frac{U_k^{\mu_2}(x)-U_k^{\mu_1}(x)}
           {(z-x)^{k+1}}\d x,
    \quad z\in\Cm.
\]
\end{lemma}

\begin{lemma}[{\cite[Lemma~2.1]{Soo26}}]
\label{lem:uniqueness}
Let $\sigma$ be a real signed Radon measure on $\R$ such that
$|\sigma|([-R,R])=o(R)$ as $R\to\infty$. Then
\[
    \int_\R\frac{\d|\sigma|(x)}{|z-x|^{k+1}}<\infty,
    \quad z\in\Cm.
\]
If
\[
    \int_\R\frac{\d\sigma(x)}{(z-x)^{k+1}}=0
    \quad\text{for every }z\in\Cm,
\]
then $\sigma=0$.
\end{lemma}

The transform identity and uniqueness result in \cite{Soo26}
are stated on $\Cp$; their versions on $\Cm$ follow by complex
conjugation.

\subsection{Multiplicative subordination}\label{subs:subord}

For $\mu\in\Pp_0(\R)$, define
\[\psi_\mu(z)
    =\int_\R\frac{zx}{1-zx}\d\mu(x),\quad\eta_\mu(z)
    =\frac{\psi_\mu(z)}{1+\psi_\mu(z)}
      =1-zF_\mu(1/z),\quad h_\mu(z)
    =\frac{\eta_\mu(z)}{z},
    \quad z\in\C\setminus\R.\]
When $\mu$ is supported on $[0,\infty)$, these functions,
as well as $G_\mu$ and $F_\mu$, extend analytically to
$\C\setminus[0,\infty)$
\cite[Proposition~2.5]{AHK26}.

Let $\mu\in\Pp_0(\R)\setminus\{\delta_0\}$ and
$\nu\in\Pp_0(\R_{\ge0})\setminus\{\delta_0\}$. By
\cite[Theorem~1.3]{AHK26}, the pair $(\mu,\nu)$ determines
unique subordination functions $\omega_\mu:\Cp\to\Cp$ and
$\omega_\nu:\Cp\to\C\setminus[0,\infty)$.
The free multiplicative convolution $\mu\boxtimes\nu$ is the
unique probability measure on $\R$ whose $\eta$-transform satisfies
\begin{equation}\label{eq:sub-eta}
    \eta_{\mu\boxtimes\nu}(z)
    =\eta_\mu(\omega_\mu(z))
    =\eta_\nu(\omega_\nu(z))
    =\frac{\omega_\mu(z)\omega_\nu(z)}{z},
    \quad z\in\Cp.
\end{equation}
The subordination functions satisfy
$\arg z\le\arg\omega_\nu(z)\le\arg z+\pi$ and
$\omega_\mu(z),\omega_\nu(z)\to0$ as $z\to0$
nontangentially. Here arguments are taken in $[0,2\pi)$. 

For $z\in\Cp$, \eqref{eq:sub-eta} gives
\begin{equation}\label{eq:fixed-point}
\begin{aligned}
    \omega_\mu(z)
    =z h_\nu(z h_\mu(\omega_\mu(z))),\quad
    \omega_\nu(z)
    =z h_\mu(\omega_\mu(z)).
\end{aligned}
\end{equation}
For each fixed $z\in\Cp$, \eqref{eq:fixed-point} shows that
$\omega_\mu(z)$ is a fixed point of the analytic self-map
$u\mapsto z h_\nu(z h_\mu(u))$ of $\Cp$.
This map is not a conformal automorphism:
when both measures are nondegenerate, this follows from
\cite[Proposition~3.3]{AHK26}; if either measure is a point
mass, the map is constant because $h_{\delta_a}=a$.
The Schwarz lemma therefore implies that $\omega_\mu(z)$
is its unique fixed point in $\Cp$.

To express \eqref{eq:sub-eta} in terms of reciprocal Cauchy
transforms, put $\lambda=\mu\boxtimes\nu$ and define
\[
    \Omega(z)=\frac1{\omega_\mu(1/z)},\quad
    \Lambda(z)=\frac1{\omega_\nu(1/z)},\quad z\in\Cm.
\]
Then $\Omega:\Cm\to\Cm$ and
$\Lambda:\Cm\to\C\setminus[0,\infty)$.
Equation~\eqref{eq:sub-eta} gives
\begin{equation}\label{eq:sub-cauchy}
\begin{aligned}
    \frac{F_\lambda(z)}{z}
    =\frac{F_\mu(\Omega(z))}{\Omega(z)}
      =\frac{F_\nu(\Lambda(z))}{\Lambda(z)},\quad
    \Omega(z)\Lambda(z)
    =\frac{z}{1-F_\lambda(z)/z},
    \quad z\in\Cm.
\end{aligned}
\end{equation}

The operation $\boxtimes$ is jointly weakly continuous on
$\Pp_0(\R)\times\Pp_0(\R_{\ge0})$
\cite[Corollary~1.2]{AHK26}.
We also record the continuity of the reciprocal subordination
functions. Suppose that $\mu_n\weak\mu$ and $\nu_n\weak\nu$,
where $\mu_n,\mu\in\Pp_0(\R)\setminus\{\delta_0\}$ and $\nu_n,\nu\in\Pp_0(\R_{\ge0})\setminus\{\delta_0\}.$
Denote the subordination functions for $(\mu_n,\nu_n)$ by
$\omega_{\mu,n},\omega_{\nu,n}$ and set
$\Omega_n(z)=1/\omega_{\mu,n}(1/z)$. Then
\begin{equation}\label{eq:subordination-continuity}
    \Omega_n(z)\longrightarrow\Omega(z),\quad z\in\Cm.
\end{equation}
Indeed, as shown in \cite[proof of Corollary~4.4]{AHK26},
$\omega_{\nu,n}(z)\to\omega_\nu(z)$ for every $z\in\Cp$,
and $\eta_{\nu_n}\to\eta_\nu$ locally uniformly on
$\C\setminus[0,\infty)$.
Since $h_{\nu_n}(z)=\eta_{\nu_n}(z)/z$, we also have
$h_{\nu_n}\to h_\nu$ locally uniformly on this domain.
As $\omega_\nu(z)\in\C\setminus[0,\infty)$, it follows that
\[
    \omega_{\mu,n}(z)
    =z h_{\nu_n}(\omega_{\nu,n}(z))
    \longrightarrow z h_\nu(\omega_\nu(z))
    =\omega_\mu(z),\quad z\in\Cp.
\]
Taking reciprocals at $1/z$ proves
\eqref{eq:subordination-continuity}.

\section{Variation along a segment}\label{sec:variation}

Fix distinct $\mu_1,\mu_2\in\Pp_0(\R)$ and
$\nu\in\Pp_0(\R_{\ge0})\setminus\{\delta_0\}$. For $0\le r\le1$,
put
\[
    \mu_r=(1-r)\mu_1+r\mu_2,\quad
    \lambda_r=\mu_r\boxtimes\nu.
\]
For $0<r<1$, the measure $\mu_r$ is nondegenerate.
We denote by $\omega_{\mu,r}$ and $\omega_{\nu,r}$ the
subordination functions for $(\mu_r,\nu)$ associated with
$\mu_r$ and $\nu$, respectively, and write $\Omega_r,\Lambda_r$
for the corresponding reciprocal subordination functions.

\begin{lemma}\label{lem:analytic}
For every $z\in\Cm$, the maps
\[
    r\longmapsto\Omega_r(z),\quad
    r\longmapsto\Lambda_r(z),\quad
    r\longmapsto G_{\lambda_r}(z)
\]
are continuously differentiable on $(0,1)$.
\end{lemma}

\begin{proof}[\sc Proof.]
Fix $r_0\in(0,1)$ and $z_0\in\Cp$, and put
$u_0=\omega_{\mu,r_0}(z_0)$. For real $r\in(0,1)$ and
$u\in\Cp$, we have
\[
    h_{\mu_r}(u)=\frac1u
    -\frac1{(1-r)G_{\mu_1}(1/u)+rG_{\mu_2}(1/u)}.
\]
Since the denominator is nonzero at $(r_0,u_0)$, the right-hand
side extends holomorphically in $(r,u)$ near $(r_0,u_0)$.
Moreover, $h_{\mu_{r_0}}(u_0)\in\Cp$ by
\cite[Proposition~3.3(1)]{AHK26}, so
$z_0h_{\mu_{r_0}}(u_0)\notin[0,\infty)$.
As $h_\nu$ is holomorphic on $\C\setminus[0,\infty)$,
\[
    \Phi(r,z,u)=z h_\nu(z h_{\mu_r}(u))
\]
is holomorphic near $(r_0,z_0,u_0)$.

By \eqref{eq:fixed-point}, $u_0=\Phi(r_0,z_0,u_0)$.
Since the self-map $u\mapsto\Phi(r_0,z_0,u)$ is not a conformal
automorphism, the Schwarz--Pick theorem at its fixed point gives
\[
    |\partial_u\Phi(r_0,z_0,u_0)|<1.
\]
Thus $1-\partial_u\Phi(r_0,z_0,u_0)\ne0$, and the holomorphic
implicit function theorem gives a local holomorphic solution
$u=u(r,z)$ to $u=\Phi(r,z,u)$, with $u(r_0,z_0)=u_0$. Since $u(r_0,z_0)=u_0\in\Cp$, continuity ensures that
$u(r,z)\in\Cp$ whenever $(r,z)$ is sufficiently close to
$(r_0,z_0)$.
For real $r\in(0,1)$ and $z\in\Cp$ in this neighborhood,
the identity $u(r,z)=\Phi(r,z,u(r,z))$ shows that
$u(r,z)$ is a fixed point in $\Cp$ of
$u\mapsto z h_\nu(z h_{\mu_r}(u))$.
By uniqueness of the fixed point in $\Cp$, we obtain
$u(r,z)=\omega_{\mu,r}(z)$.

The relation
$\omega_{\nu,r}(z)=z h_{\mu_r}(\omega_{\mu,r}(z))$
gives the same local holomorphic dependence for $\omega_{\nu,r}$.
Taking reciprocals at $1/z$ therefore proves the assertion for
$\Omega_r$ and $\Lambda_r$. Finally, \eqref{eq:sub-cauchy} gives
\[
    G_{\lambda_r}(z)
    =\frac{\Omega_r(z)}{z}
      \left(
      (1-r)G_{\mu_1}(\Omega_r(z))
      +rG_{\mu_2}(\Omega_r(z))
      \right),
\]
whose right-hand side is continuously differentiable in $r$.
\end{proof}

\begin{proposition}\label{prop:variation}
For $0<r<1$ and $z\in\Cm$,
\begin{equation}\label{eq:variation}
    \partial_rG_{\lambda_r}(z)
    =\Omega_r'(z)
      \left(G_{\mu_2}(\Omega_r(z))-G_{\mu_1}(\Omega_r(z))\right).
\end{equation}
\end{proposition}

\begin{proof}[\sc Proof.]
By Lemma~\ref{lem:analytic}, we may differentiate the subordination
identities with respect to $r$.
Put $P_r(z)=zG_{\lambda_r}(z)$. By \eqref{eq:sub-cauchy},
\[
    P_r(z)
    =\Omega_r(z)G_{\mu_r}(\Omega_r(z))
      =\Lambda_r(z)G_\nu(\Lambda_r(z)),\quad
    \Omega_r(z)\Lambda_r(z)
    =\frac{zP_r(z)}{P_r(z)-1}.
\]
We suppress the argument $z$ below and use primes for
differentiation in $z$.

Since $G_{\mu_r}=(1-r)G_{\mu_1}+rG_{\mu_2}$, differentiating
the identities for $P_r$ with respect to $r$ and $z$ gives
\[
    (\partial_rP_r)\Omega_r'-P_r'\partial_r\Omega_r
    =\Omega_r\Omega_r'
      \left(G_{\mu_2}(\Omega_r)-G_{\mu_1}(\Omega_r)\right),\quad
    (\partial_rP_r)\Lambda_r'-P_r'\partial_r\Lambda_r=0.
\]
Multiplying these identities by $\Lambda_r$ and $\Omega_r$,
respectively, and adding, we obtain
\[
    (\partial_rP_r)\partial_z(\Omega_r\Lambda_r)
      -P_r'\partial_r(\Omega_r\Lambda_r)=\Omega_r\Lambda_r\Omega_r'
      \left(G_{\mu_2}(\Omega_r)-G_{\mu_1}(\Omega_r)\right).
\]
Substituting $\Omega_r\Lambda_r=zP_r/(P_r-1)$, the terms
containing $P_r'\partial_rP_r$ cancel, leaving
\[
    \frac{P_r}{P_r-1}\partial_rP_r
    =\frac{zP_r}{P_r-1}\Omega_r'
      \left(G_{\mu_2}(\Omega_r)-G_{\mu_1}(\Omega_r)\right).
\]
Cancelling the nonzero factor $P_r/(P_r-1)$ gives
\[
    \partial_rP_r
    =z\Omega_r'
      \left(G_{\mu_2}(\Omega_r)-G_{\mu_1}(\Omega_r)\right).
\]
Since $G_{\lambda_r}(z)=P_r(z)/z$, this proves
\eqref{eq:variation}.
\end{proof}
For the kernel construction below, we express this variation
in terms of the potentials. If
$\mu_1,\mu_2\in\Pp_{k-1}(\R)$ have matching moments through
order $k-1$, then \eqref{eq:variation} and
Lemma~\ref{lem:potentials} give
\begin{equation}\label{eq:variation-potentials}
    \partial_rG_{\lambda_r}(z)
    =k!\int_\R
      \left(U_k^{\mu_2}(t)-U_k^{\mu_1}(t)\right)
      \frac{\Omega_r'(z)}{(\Omega_r(z)-t)^{k+1}}\d t,
    \quad 0<r<1,\quad z\in\Cm.
\end{equation}

\section{Positive subordination kernels}\label{sec:kernels}
In this section, let $\mu\in\Pp_0(\R)\setminus\{\delta_0\}$ and
$\nu\in\Pp_1(\R_{\ge0})\setminus\{\delta_0\}$.
We construct finite positive measures representing the kernel
in \eqref{eq:variation-potentials} and prove an integrated
bound depending only on $m_\nu$.

Since $\nu$ has finite first moment, dominated convergence gives
$h_\nu(z)\to m_\nu$ as $z\to0$ with $\arg z$ bounded away from
$0$ and $2\pi$. The subordination limits and the angular bound
in \eqref{eq:sub-eta} therefore yield
\[
    \frac{\Omega(z)}{z}
    =\frac1{h_\nu(\omega_\nu(1/z))}
    \longrightarrow\frac1{m_\nu}
\]
as $z\to\infty$ nontangentially in $\Cm$.

Since $m_\nu>0$, the reciprocal Cauchy transform characterization
\cite[Proposition~2.3]{AHK26}, reflected across $\R$, identifies
$m_\nu(\Omega-t)$ as the reciprocal Cauchy transform of a
probability measure for every $t\in\R$. Multiplying this measure
by $m_\nu$, we obtain a unique finite positive measure $\rho_t$
satisfying
\begin{equation}\label{eq:rho}
    G_{\rho_t}(z)=\frac1{\Omega(z)-t},
    \quad z\in\Cm,\quad \rho_t(\R)=m_\nu.
\end{equation}

For $t\in\R$ and an integer $j\ge1$, let
$X_1,\ldots,X_j$ be independent random variables with common
law $m_\nu^{-1}\rho_t$. Independently, let
$(\Theta_1,\ldots,\Theta_j)$ have the
$\operatorname{Dirichlet}(1,\ldots,1)$ distribution, with
$\Theta_1=1$ when $j=1$. Define
\begin{equation}\label{eq:kernel-definition}
    K_t^{(j)}
    =m_\nu^j\Law\left(\sum_{i=1}^j\Theta_iX_i\right).
\end{equation}
These are finite positive measures satisfying
\begin{equation}\label{eq:kernel-mass}
    K_t^{(j)}(\R)=m_\nu^j,\quad K_t^{(1)}=\rho_t.
\end{equation}

For fixed $x_1,\ldots,x_j\in\R$, the Hermite--Genocchi formula
\cite[Eq.~(A.3)]{Les17}, with normalized Dirichlet weights and
$g(x)=(z-x)^{-1}$, gives
\[
    \prod_{i=1}^j\frac1{z-x_i}
    =\mathbb E_\Theta
      \left[
          \frac1{(z-\sum_{i=1}^j\Theta_i x_i)^j}
      \right],
    \quad z\in\Cm.
\]
Applying this identity at $x_i=X_i$ and taking expectations,
independence and \eqref{eq:kernel-definition} yield
\begin{equation}\label{eq:kernel-power}
\begin{aligned}
    \int_\R\frac{\d K_t^{(j)}(x)}{(z-x)^j}
    =m_\nu^j\mathbb E
      \left[\prod_{i=1}^j\frac1{z-X_i}\right]=G_{\rho_t}(z)^j
      =\frac1{(\Omega(z)-t)^j}.
\end{aligned}
\end{equation}
The preceding expectations are finite since
$|z-x|\ge|\Im z|$ for $x\in\R$.
Moreover, the integrand in \eqref{eq:kernel-power} and its
$z$-derivative are bounded uniformly in $x$ when $z$ ranges
over a compact subset of $\Cm$.
Since $K_t^{(j)}$ is finite, we may differentiate
\eqref{eq:kernel-power} under the integral.
Taking $j=k$ and cancelling the factor $-k$ gives
\begin{equation}\label{eq:kernel-rep}
    \int_\R\frac{\d K_t^{(k)}(x)}{(z-x)^{k+1}}
    =\frac{\Omega'(z)}{(\Omega(z)-t)^{k+1}},
    \quad z\in\Cm.
\end{equation}

\begin{lemma}\label{lem:kernel-continuity}
Suppose that $(\mu_n,\nu_n)$ and $(\mu,\nu)$ satisfy the
assumptions of this section, with $\mu_n\weak\mu$,
$\nu_n\weak\nu$ and $m_{\nu_n}\to m_\nu$.
Write $K_{n,t}^{(k)}$ for the kernel associated with
$(\mu_n,\nu_n)$. If $t_n\to t$, then
\[
    K_{n,t_n}^{(k)}\weak K_t^{(k)}.
\]
\end{lemma}

\begin{proof}[\sc Proof.]
Let $\Omega_n$ and $\rho_{n,t}$ be associated with
$(\mu_n,\nu_n)$.
By \eqref{eq:rho} and \eqref{eq:subordination-continuity},
\[
    G_{\rho_{n,t_n}}(z)=\frac1{\Omega_n(z)-t_n}
    \longrightarrow\frac1{\Omega(z)-t}=G_{\rho_t}(z).
\]
Since the total masses converge to $m_\nu>0$, the continuity
theorem for Cauchy transforms \cite[Theorem~2.2]{AHK26},
applied to the normalized measures, gives
$\rho_{n,t_n}\weak\rho_t$.
The construction in \eqref{eq:kernel-definition} then gives
$K_{n,t_n}^{(k)}\weak K_t^{(k)}$ by convergence of product
measures, the continuous mapping theorem, and
$m_{\nu_n}^k\to m_\nu^k$.
\end{proof}
In particular, $t\mapsto K_t^{(k)}$ is weakly continuous.
Along the segment of Section~\ref{sec:variation}, write
$K_{r,t}^{(k)}$ for the corresponding kernels.
Since $r\mapsto\mu_r$ is weakly continuous and $\nu$ is fixed,
Lemma~\ref{lem:kernel-continuity} shows that, for every bounded
continuous function $f:\R\to\C$, the map
\[
    (r,t)\longmapsto\int_\R f(x)\d K_{r,t}^{(k)}(x)
\]
is continuous on $(0,1)\times\R$.
A monotone class argument therefore shows that
$t\mapsto K_t^{(k)}(A)$ and
$(r,t)\mapsto K_{r,t}^{(k)}(A)$ are Borel measurable
for every Borel set $A\subset\R$.

We next establish the integrated bound needed in the main proof.

\begin{lemma}\label{lem:kernel-bound}
For every nonnegative Borel function $\varphi$ on $\R$,
\[
    \int_\R\left(\int_\R\varphi(x)\d K_t^{(k)}(x)\right)\d t
    \le m_\nu^{k-1}\int_\R\varphi(x)\d x.
\]
\end{lemma}

\begin{proof}[\sc Proof.]
Fix $R>0$. For every Borel set $A\subset\R$, define
\[
    N_R(A)=\int_{-R}^R K_t^{(k)}(A)\d t.
\]
The measurability established above and \eqref{eq:kernel-mass}
show that $N_R$ is a finite positive Borel measure.
Define the bounded Borel function
\[
    f_R(x)=K_R^{(k-1)}((x,\infty))
           -K_{-R}^{(k-1)}((x,\infty)).
\]
Both kernels have mass $m_\nu^{k-1}$, so
$|f_R|\le m_\nu^{k-1}$ and $f_R(x)\to0$ as $x\to\pm\infty$.

For $z\in\Cm$, Fubini's theorem and \eqref{eq:kernel-power} give
\[
\begin{aligned}
    \int_\R\frac{f_R(x)}{(z-x)^k}\d x
    &=\int_\R
      \left(\int_{-\infty}^x\frac{\d y}{(z-y)^k}\right)
      \left(\d K_R^{(k-1)}(x)-\d K_{-R}^{(k-1)}(x)\right)\\
    &=\frac1{k-1}\int_\R
      \frac{\d K_R^{(k-1)}(x)-\d K_{-R}^{(k-1)}(x)}
           {(z-x)^{k-1}}\\
    &=\frac{(\Omega(z)-R)^{1-k}
            -(\Omega(z)+R)^{1-k}}{k-1}\\
    &=\int_{-R}^R\frac{\d t}{(\Omega(z)-t)^k}\\
    &=\int_{-R}^R
      \left(\int_\R\frac{\d K_t^{(k)}(x)}{(z-x)^k}\right)\d t\\
    &=\int_\R\frac{\d N_R(x)}{(z-x)^k}.
\end{aligned}
\]
Since $N_R$ is finite and $f_R$ vanishes at both infinities,
the signed Radon measure $N_R-f_R(x)\d x$ has total variation
$o(L)$ on $[-L,L]$ as $L\to\infty$.
Differentiating the preceding identity once in $z$ and applying
Lemma~\ref{lem:uniqueness}, we obtain $ \d N_R(x)=f_R(x)\d x.$
Consequently, $N_R\le m_\nu^{k-1}\d x$.
Integrating $\varphi$ and letting $R\to\infty$ proves the
assertion by monotone convergence.
\end{proof}

\section{Proof of Theorem~\ref{thm:main}}\label{sec:proof}

We now prove Theorem~\ref{thm:main}. The cases $\mu_1=\mu_2$ and
$\nu=\delta_0$ are immediate, so we assume $\mu_1\ne\mu_2$ and
$\nu\ne\delta_0$. We use the segment notation of Section~\ref{sec:variation},
so that $\lambda_0=\mu_1\boxtimes\nu$ and
$\lambda_1=\mu_2\boxtimes\nu$.
By Proposition~\ref{prop:DLS}, $\mu_1$ and $\mu_2$ have matching
moments through order $k-1$, and $U_k^{\mu_1}(t)\le U_k^{\mu_2}(t)$ for all $t\in\R$.

\medskip
\noindent\textbf{Claim} $(\star)$.
$U_k^{\lambda_0}(t)\le U_k^{\lambda_1}(t)$ for every $t\in\R$.

\begin{proof}[\sc Proof of Claim~$(\star)$]
For every Borel set $A\subset\R$, define
\begin{equation}\label{eq:Mdef}
    M_k(A)=\int_0^1\int_\R
    \left(U_k^{\mu_2}(t)-U_k^{\mu_1}(t)\right)
    K_{r,t}^{(k)}(A)\d t\d r.
\end{equation}
The joint measurability established after
Lemma~\ref{lem:kernel-continuity} and the nonnegativity of
$U_k^{\mu_2}-U_k^{\mu_1}$ ensure that \eqref{eq:Mdef}
defines a positive Borel measure.

Lemma~\ref{lem:potentials} gives
$U_k^{\mu_2}-U_k^{\mu_1}\in C_0(\R)$. Given $\varepsilon>0$, choose $T>0$ such that
\[
    0\le U_k^{\mu_2}(t)-U_k^{\mu_1}(t)<\varepsilon
    \quad\text{whenever }|t|>T.
\]
We split the integral in \eqref{eq:Mdef} into the regions
$|t|\le T$ and $|t|>T$. The mass identity
\eqref{eq:kernel-mass} controls the first contribution, while Lemma~\ref{lem:kernel-bound} controls the second. Thus, for $R>0$,
\[
\begin{aligned}
    M_k([-R,R])
    &\le
       2T m_\nu^k
       \left\|U_k^{\mu_2}-U_k^{\mu_1}\right\|_\infty
       +\varepsilon\int_0^1\int_\R
        K_{r,t}^{(k)}([-R,R])\d t\d r\\
    &\le
       2T m_\nu^k
       \left\|U_k^{\mu_2}-U_k^{\mu_1}\right\|_\infty
       +2\varepsilon R m_\nu^{k-1}.
\end{aligned}
\]
It follows that $M_k$ is finite on compact sets and so is a
Radon measure. Dividing by $R$ and letting $R\to\infty$ gives
\[
    \limsup_{R\to\infty}\frac{M_k([-R,R])}{R}
    \le2\varepsilon m_\nu^{k-1}.
\]
Since $\varepsilon>0$ is arbitrary, we deduce that $M_k([-R,R])=o(R)$ as $R\to\infty$.

Fix $z\in\Cm$. Lemma~\ref{lem:uniqueness} gives
$\int_\R|z-x|^{-k-1}\d M_k(x)<\infty$,
which justifies Fubini's theorem below:
\begin{equation}\label{eq:Mtransform}
\begin{aligned}
    k!\int_\R\frac{\d M_k(x)}{(z-x)^{k+1}}
    &=
    k!\int_0^1\int_\R
    \left(U_k^{\mu_2}(t)-U_k^{\mu_1}(t)\right)
    \left(
        \int_\R\frac{\d K_{r,t}^{(k)}(x)}{(z-x)^{k+1}}
    \right)\d t\d r\\
    &=
    k!\int_0^1\int_\R
    \left(U_k^{\mu_2}(t)-U_k^{\mu_1}(t)\right)
    \frac{\Omega_r'(z)}{(\Omega_r(z)-t)^{k+1}}
    \d t\d r\quad(\text{by \eqref{eq:kernel-rep}})\\
    &=
    \int_0^1\partial_rG_{\lambda_r}(z)\d r
    \quad\text{(by \eqref{eq:variation-potentials})}\\
    &=
    G_{\lambda_1}(z)-G_{\lambda_0}(z).
\end{aligned}
\end{equation}
The last equality follows by integrating over
$[\varepsilon,1-\varepsilon]$ and letting $\varepsilon\downarrow0$,
using the preceding absolute-integrability estimate and the weak continuity of $r\mapsto\lambda_r$.

To apply Lemma~\ref{lem:potentials} to $\lambda_0,\lambda_1$,
we first verify that they belong to $\Pp_p(\R)$, where $p=k-1$.
For bounded freely independent operators $X=X^*$ and $Y\ge0$
in a finite von Neumann algebra with a faithful normal tracial
state $\tau$,
write $\|S\|_q=\tau(|S|^q)^{1/q}$.
Using the polar decomposition of $X$, we obtain
\begin{equation}\label{eq:moment-estimate}
\begin{aligned}
    \|Y^{1/2}XY^{1/2}\|_p^p
    &\le\||X|^{1/2}Y^{1/2}\|_{2p}^{2p}
    \quad\text{(H\"older, \cite[Theorem~4.2(i)]{FK86})}\\
    &\le\tau(|X|^pY^p)
    \quad\text{(Araki--Lieb--Thirring,
    \cite[Corollary~3]{Kos92})}\\
    &=\|X\|_p^p\|Y\|_p^p
    \quad\text{(freeness)}.
\end{aligned}
\end{equation}
Applying this estimate to bounded spectral truncations, the
weak continuity of $\boxtimes$ \cite[Corollary~1.2]{AHK26}
and the portmanteau theorem extend it to freely independent
affiliated operators $X=X^*$ and $Y\ge0$ in $L^p$.
In particular, $\lambda_0,\lambda_1\in\Pp_{k-1}(\R)$.

For compactly supported measures, the moment--cumulant formula
and the product formula for free cumulants
\cite[Proposition~11.4 and Theorem~14.4]{NS06}
express $m_j(\mu_i\boxtimes\nu)$ as a universal polynomial in
the moments of $\mu_i$ and $\nu$ through order $j$.
For $1\le j\le p$, these identities extend by successive spectral
truncation: H\"older's inequality controls truncation of $Y$ for
bounded $X$, and \eqref{eq:moment-estimate} controls truncation
of $X$ with $Y$ fixed. The resulting $L^p$ convergence implies
convergence of moments through order $p$, again by H\"older's
inequality; the input moments converge by dominated convergence.
Since $\mu_1$ and $\mu_2$ have matching moments through
order $k-1$, we obtain
\begin{equation}\label{eq:output-moments}
    m_j(\lambda_0)=m_j(\lambda_1),
    \quad 1\le j\le k-1.
\end{equation}

Lemma~\ref{lem:potentials} now gives
$U_k^{\lambda_1}-U_k^{\lambda_0}\in C_0(\R)$ and
\begin{equation}\label{eq:outputpotentialtransform}
    G_{\lambda_1}(z)-G_{\lambda_0}(z)
    =
    k!\int_\R
    \frac{U_k^{\lambda_1}(x)-U_k^{\lambda_0}(x)}
         {(z-x)^{k+1}}\d x,
    \quad z\in\Cm.
\end{equation}
Comparing \eqref{eq:Mtransform} and
\eqref{eq:outputpotentialtransform}, the real signed Radon measure
\[
    \d\sigma(x)=\d M_k(x)
    -\left(U_k^{\lambda_1}(x)-U_k^{\lambda_0}(x)\right)\d x
\]
satisfies
\[
    \int_\R\frac{\d\sigma(x)}{(z-x)^{k+1}}=0,
    \quad z\in\Cm.
\]
Since $M_k([-R,R])=o(R)$ and
$U_k^{\lambda_1}-U_k^{\lambda_0}\in C_0(\R)$,
we also have $|\sigma|([-R,R])=o(R)$.
Lemma~\ref{lem:uniqueness} therefore gives $\sigma=0$, so
\[
    \d M_k(x)
    =\left(U_k^{\lambda_1}(x)-U_k^{\lambda_0}(x)\right)\d x.
\]
The positivity of $M_k$ implies
$U_k^{\lambda_0}\le U_k^{\lambda_1}$ almost everywhere.
By continuity, the inequality holds everywhere,
proving Claim~$(\star)$.
\end{proof}

By Claim~$(\star)$, together with \eqref{eq:output-moments},
Proposition~\ref{prop:DLS} yields
$\mu_1\boxtimes\nu\kcx\mu_2\boxtimes\nu$,
as required.\qed


\begin{thebibliography}{99}

\bibitem[AHK26]{AHK26}
O.~Arizmendi, T.~Hasebe, and Y.~Kitagawa,
\newblock Free multiplicative convolution with an arbitrary measure
on the real line,
\newblock \emph{Probab. Theory Related Fields} \textbf{196} (2026),
1761--1813,
doi:10.1007/s00440-026-01483-z.

\bibitem[BV93]{BV93}
H.~Bercovici and D.~Voiculescu,
\newblock Free convolution of measures with unbounded support,
\newblock \emph{Indiana Univ. Math. J.} \textbf{42} (1993), no.~3,
733--773.

\bibitem[Boc91]{Boca91}
F.~Boca,
\newblock Free products of completely positive maps and spectral sets,
\newblock \emph{J. Funct. Anal.} \textbf{97} (1991), no.~2,
251--263.

\bibitem[DLS98]{DLS98}
M.~Denuit, C.~Lef\`evre, and M.~Shaked,
\newblock The $s$-convex orders among real random variables,
with applications,
\newblock \emph{Math. Inequal. Appl.} \textbf{1} (1998), no.~4,
585--613.

\bibitem[FK86]{FK86}
T.~Fack and H.~Kosaki,
\newblock Generalized $s$-numbers of $\tau$-measurable operators,
\newblock \emph{Pacific J. Math.} \textbf{123} (1986), no.~2,
269--300.

\bibitem[Hei26]{Hei26}
O.~Hein\"avaara,
\newblock Convolution comparison measures,
\newblock arXiv:2602.10373, 2026.

\bibitem[Hia87]{Hiai87}
F.~Hiai,
\newblock Majorization and stochastic maps in von Neumann algebras,
\newblock \emph{J. Math. Anal. Appl.} \textbf{127} (1987), no.~1,
18--48.

\bibitem[Kos92]{Kos92}
H.~Kosaki,
\newblock An inequality of Araki--Lieb--Thirring
(von Neumann algebra case),
\newblock \emph{Proc. Amer. Math. Soc.} \textbf{114} (1992), no.~2,
477--481.

\bibitem[Les17]{Les17}
M.~Lesch,
\newblock Divided differences in noncommutative geometry:
rearrangement lemma, functional calculus and expansional formula,
\newblock \emph{J. Noncommut. Geom.} \textbf{11} (2017), no.~1,
193--223.

\bibitem[NS06]{NS06}
A.~Nica and R.~Speicher,
\newblock \emph{Lectures on the Combinatorics of Free Probability},
\newblock London Mathematical Society Lecture Note Series, vol.~335,
Cambridge University Press, Cambridge, 2006.

\bibitem[Soo26]{Soo26}
C.-X.~Soong,
\newblock Higher-order convex orders under free additive convolution,
\newblock arXiv:2609.35901, 2026.

\bibitem[Str65]{Strassen65}
V.~Strassen,
\newblock The existence of probability measures with given marginals,
\newblock \emph{Ann. Math. Statist.} \textbf{36} (1965), no.~2,
423--439.

\bibitem[Voi87]{Voi87}
D.~Voiculescu,
\newblock Multiplication of certain non-commuting random variables,
\newblock \emph{J. Operator Theory} \textbf{18} (1987), no.~2,
223--235.

\end{thebibliography}
\end{document}